\documentclass[12pt, reqno]{amsart}
\usepackage{amsmath, amsthm, amscd, amsfonts, amssymb, graphicx, color}
\usepackage[bookmarksnumbered, colorlinks, plainpages]{hyperref}

\newtheorem{thm}{Theorem}[section]
 \newtheorem{cor}[thm]{Corollary}

 \theoremstyle{definition}
 
 \theoremstyle{remark}
 
 \newtheorem*{ex}{Example}
 \numberwithin{equation}{section}
\begin{document}
\title[A note on Bernstein property ]{A note on Bernstein property of a fourth order complex partial differential equations}
\author{Sa\"{\i}d Asserda}
\address{Ibn tofail university , faculty of sciences, department of mathematics, PO 242 Kenitra Morroco}
\email{asserda-said@univ-ibntofail.ac.ma}
\subjclass[2010]{Primary 58E20, Secondary 53C21}
\keywords{Ricci curvature, Bernstein property}
\date{\today}
\begin{abstract} For a smooth strictly
plurisubharmonic function $u$ on a open set
$\Omega\subset\mathbb{C}^{n}$ and $F$ a $C^{1}$ nondecreasing
function on $\mathbf{R}^{*}_{+}$, we investigate the complex partial
differential equations
$$\Delta_{g}\log\det(u_{i\bar j})=F(\det(u_{i\bar
j}))\Vert\nabla_{g}\log\det(u_{i\bar j})\Vert_{g}^{2},$$ where
$\Delta_{g}$, $\Vert . \Vert_{g}$ and $\nabla_{g}$ are the
Laplacian, tensor norm and the Levi-Civita connexion , respectively,
with respect to the K\"ahler metric $g=\partial\bar\partial u$. We
show that the above PDE's has a Bernstein property, i.e
$\det(u_{i\bar j})=\hbox{constant}$ on $\Omega$, provided that $g$
is complete, the Ricci curvature of $g$ is bounded below and $F$
satisfies
 $\inf_{t\in\mathbf{R}^{+}}(2tF^{'}(t)+{F(t)^{2}\over n})>{1\over 4}$
and $F(\max_{B(R)}\det u_{i\bar j})=o(R).$
\end{abstract}
\maketitle
\section{Introduction}
In this  note, we investigate a class of fourth order complex
nonlinear partial differential equations for a strictly
plurisubharmonic function $u$ on a open set
$\Omega\subset\mathbb{C}^{n}$ :
\begin{equation}\Delta_{g}\log\det(u_{i\bar j})=F(\det(u_{i\bar j}))\Vert\nabla_{g}\log\det(u_{i\bar j})\Vert_{g}^{2}
\end{equation}
where $\Delta_{g}$, $\Vert . \Vert_{g}$ and $\nabla_{g}$ are the
Laplacian, tensor norm and the Levi-Civita connexion , respectively,
with respect to the K\"ahler metric $g=\partial\bar\partial u$. For
$n\geq 2$, the equation (1.1) is the Euler-Lagrange equation of the
variational problem of the functional
\begin{equation}{\mathcal V}_{\Phi}(u)=\int_{\Omega}\Phi(\det(u_{i\bar j}))dV_{e}
\end{equation}
 for some $\Phi\in C^{4}(\mathbf{R}^{*}_{+})$. If  $\Phi(t)=t^{1-\beta}, \ \beta\not=1$, the Euler-Lagrange equation of ${\mathcal V}_{\Phi}$ is
\begin{equation}
\Delta_{g}\log\det(u_{i\bar
j})=\beta\Vert\nabla_{g}\log\det(u_{i\bar j})\Vert_{g}^{2}
\end{equation}
In [1], Chen and Li showed the following  Bernstein property  of the
solution of (3)
\begin{thm} Let
$\Omega\subset\mathbb{C}^{n}$ and $u\in C^{2}(\Omega)$ a strictly
plurisubharmonic solution of (3).
 Assume that \\
(i) $(\Omega,g)$ is complete,\\
(ii) $Ric(g)\geq -C_{1}\ ( K\geq 0)$,\\
(ii) $\det(u_{i\bar j})< C_{2}$\\
where $C_{1},C_{2} >0$ are real constants. If
$\beta>{\sqrt{9n^{2}+4n}+3n\over 4}$ then $$ \det(u_{i\bar
j})=\hbox{const}\ \hbox{on}\ \Omega.
$$
\end{thm} 
For more details
about Bernstein property, one can see the paper [1] and references
therein. In this note, we will show a Bernstein property of the
solutions of (1.1). Our main result is the following.
\begin{thm} Let $\Omega\subset\mathbb{C}^{n}$
be an open set and $u\in C^{2}(\Omega)$ a strictly plurisubharmonic
solution of (1).
 Assume that \\
(i) $(\Omega,g)$ is complete,\\
(ii) $Ric(g)\geq -K\ ( K\geq 0)$,\\
(ii) $F\in C^{1}(]0,+\infty[)$ nondecreasing satisfying
\begin{equation}\label{testequation}
\inf_{t>0}\left(2tF^{'}(t)+{F(t)^{2}\over n}\right)>{1\over
4}\quad\hbox{and}\quad F(\max_{B(R)}\det(u_{i\bar j}))=o(R)
\end{equation}
where $B(R)=\{ z\in \Omega, d_{g}(z,z_{0})<R\}$. Then $\det(u_{i\bar
j})$ is constant on $\Omega$.
\end{thm} 
If $F(t)=\beta$,  the
following corollary  extends theorem (1.1).
\begin{cor} Let $\Omega\subset\mathbb{C}^{n}$ be an open set
 and $u\in C^{2}(\Omega)$ a strictly plurisubharmonic
solution of (3). Assume that  $(\Omega,g)$ is complete and
$Ric(g)\geq -K\ ( K\geq 0)$. If $\vert\beta\vert>{\sqrt{n}\over 2}$
then $\det(u_{i\bar j})$ is
constant on $\Omega$.
\end{cor}
\begin{ex}
Example of function $F$ which satisfies  conditions of Theorem 1.2 : \\
1- $F(t)=\log^{2\alpha}t$ where $\alpha\in]\alpha_{0}(n),{1\over 2}],\ \ \alpha_{0}(n)>0$.\\
2- $F(t)=(\log t)\log(\alpha+(\log t)^{2})$ where $\alpha> e^{1\over 4}$.\\
3- $F(t)=f(\log t)$ where $f$ is a solution of the Riccati equation
$2y^{'}+{y^{2}\over n}=\alpha,\ \alpha>{1\over 4}$.
\end{ex}
\section{Proof of theorem 1.2}
\begin{proof}
We denote $\Psi=\vert\nabla U\vert^{2}$ where $U=\log\det(u_{i\bar
j})$ and $f(t)=F(e^{t})$. By  Bochner formulae\\
\begin{eqnarray*}
\Delta_{g}\Psi&=&2f^{'}(U)\Psi^{2}+f(U)(<\nabla_{g}\Psi,\nabla_{g}
U>+<\nabla_{g} U,\nabla_{g}\Psi>)\\ &{}&+\sum u^{i\bar j}u^{k\bar
l}U_{i\bar l}U_{\bar j\bar l} +\sum u^{i\bar j}u^{k\bar l}U_{i\bar
l}U_{{\bar j}k}+\sum u^{i\bar j}u^{m\bar s}U_{i\bar l}U_{m}U_{\bar
j}R_{i\bar s}\qquad(2.1)
\end{eqnarray*}
where $(u^{i\bar j})$ is the inverse matrix of $(u_{i\bar j})$ and
$R_{i\bar s}$ is the Ricci curvature of $g$ 
$$ 
R_{i\bar
s}=-{\partial ^{2}\log\det(u_{k\bar l})\over\partial
z_{i}\partial\bar z_{s}}={\partial^{2}\log U\over\partial
z_{i}\partial\bar z_{s}}
$$ 
Let $z\in \Omega$ fixed, we can choose a
normal coordinates $(z_{1},\cdots,z_{n})$ at $z$ such that
$$u_{i\bar j}=\delta_{ij},\ \ U_{1}=U_{\bar 1},\ \ U_{i}=U_{\bar
i}=0,\ \ \forall i>1$$ At $z$ the formulae (2.1) become\\
\begin{eqnarray*}
\qquad{\Delta_{g}\Psi\over\Psi}=2f^{'}(U)\Psi &{}&+(f(U)-{1\over 2})({\Psi_{1}\over\Psi}U_{\bar 1}+U_{1}{\Psi_{\bar 1}\over\Psi})+{1\over 2}({\Psi_{1}\over\Psi}U_{\bar 1}+U_{1}{\Psi_{\bar 1}\over\Psi})\\
&{}&+{\sum U_{ik}U_{\bar j\bar k}\over\Psi}+{\sum U_{i\bar k}U_{{\bar
i}k}\over\Psi} -{U_{1\bar 1}U_{1}U_{\bar
1}\over\Psi}\qquad\qquad(2.2)
\end{eqnarray*}\\
Since $\sum U_{ik}U_{\bar i\bar k}\geq U_{11}U_{\bar 1\bar 1}$ and\\
\begin{eqnarray*}
{1\over 2}({\Psi_{1}\over\Psi}U_{\bar 1}+U_{1}{\Psi_{\bar 1}\over\Psi})-{U_{1\bar 1}U_{1}U_{\bar 1}\over\Psi}&{}&={1\over 2\Psi}2\Re(U_{11}U_{\bar 1}U_{\bar 1}\\
\geq{U_{11}U_{\bar 1\bar 1}\over\Psi}-{1\over 4}\Psi
\end{eqnarray*}\\
we obtain\\
\begin{equation*}
{1\over 2}({\Psi_{1}\over\Psi}U_{\bar 1}+U_{1}{\Psi_{\bar
1}\over\Psi})-{U_{1\bar 1}U_{1}U_{\bar 1}\over\Psi}+{\sum
U_{ik}U_{\bar i\bar k}\over\Psi}\geq-{1\over
4}\Psi\qquad\qquad\qquad\quad(2.3)
\end{equation*}\\
Let $G$ be the function defined on $B(R)=B_{g}(z_{0},R)$ by
$$ G(z)=(R^{2}-r(z)^{2})^{2}\Psi$$
where $r(z)=d_{g}(z,z_{0})$. Since $G$ vanishes at the boundary
$\partial B(R)$, the maximum of $G$ is attained on $w\in B(R)$. By
maximum principle\\
\begin{eqnarray*}
\nabla_{g} G(w)&=&0\\
\Delta_{g} G(w)&\leq &0
\end{eqnarray*}\\
A computation at $w$ gives\\
\begin{eqnarray*}
\qquad{\Psi_{i}\over\Psi}&=&{4rr_{i}\over R^{2}-r^{2}}\qquad\qquad\qquad\qquad\qquad\qquad\qquad\qquad\qquad\qquad(2.4)\\
\qquad\qquad\Delta_{g} G &=& 2\Psi\Vert\nabla_{g}(R^{2}-r^{2})\Vert^{2}+2(R^{2}-r^{2})\Psi\Delta_{g}(R^{2}-r^{2})\\
&{}&+2(R^{2}-r^{2})\sum(R^{2}-r^{2})_{i}\Psi_{i}
+(R^{2}-r^{2})^{2}\Delta_{g}\Psi\\
&=& 2\Psi r^{2}\Vert\nabla_{g}r\Vert^{2}-2\Psi(R^{2}-r^{2})\Delta_{g} r^{2}+(R^{2}-r^{2})^{2}\Delta_{g}\Psi\\
&{}& \quad-4r(R^{2}-r^{2})\sum_{i}r_{i}\Psi_{i}\leq
0\qquad\qquad\qquad\qquad\qquad\qquad(2.5)
\end{eqnarray*}
Since $\Vert\nabla_{g} r\Vert_{g}=1$, (2.4) and (2.5) implies
\begin{equation*}
\qquad\qquad\qquad\qquad\quad{\Delta_{g}\Psi\over\Psi}\leq{2\Delta_{g}
r^{2}\over
R^{2}-r^{2}}+{14r^{2}\over(R^{2}-r^{2})^{2}}\qquad\qquad\qquad\qquad\qquad(2.6)
\end{equation*}
Let $\eta >$ fixed such that $\inf_{t\in\mathbb
R}(2f^{'}(t)+{f(t)^{2}\over n})>2\eta+{1\over 4}$. By Shwarz
inequality
\begin{eqnarray*}
\qquad(f(U)-{1\over 2})({\Psi_{1}\over\Psi}U_{\bar 1}+U_{1}{\Psi_{\bar 1}\over\Psi})&=&2\Re(U_{1}{4rr_{1}(f(U)-{1\over 2})\over R^{2}-r^{2}})\\
&\geq&-\eta\Psi-{16(f(U)-{1\over
2})^{2}r^{2}\over\eta(^{2}-r^{2})^{2}}\quad\qquad\qquad(2.7)
\end{eqnarray*}
Since $\Delta U=f(U)$ we have
\begin{equation*}
\qquad\quad\quad{\sum U_{i\bar k}U_{{\bar
i}k}\over\Psi}\geq{\sum(U_{i\bar i})^{2}\over\Psi}\geq{(\sum
U_{i\bar i})^{2}\over n\Psi}={(\Delta U)^{2}\over
n\Psi}={f(U))^{2}\over n}\Psi\qquad\quad(2.8)
\end{equation*}
The inequalities (2.2), (2.7) and (2.8) give
\begin{equation*}
\qquad\qquad\qquad\qquad\qquad{\Delta_{g}\Psi\over\Psi}\geq\eta\Psi-{16(f(U)-{1\over
2})^{2}r^{2}\over\eta(R^{2}-r^{2})^{2}}\qquad\qquad\qquad\qquad(2.9)
\end{equation*}
By Laplacian comparison theorem [3] : $\Delta_{g} r^{2}\leq
2n+2(n-1)\sqrt{K}r$. Since $\vert f(U)-{1\over 2}\vert^{2}\leq
2f(U)^{2}+{1\over 2}$ and $f$ is nondecreasing, the inequalities
(2.6) and (2.9) give on $B(R)$
\begin{equation*}
(R^{2}-r^{2})^{2}\Psi\leq G(w)\leq
C_{1}(1+R)R^{2}+C_{2}(f(\max_{B(R)}U))^{2}R^{2}
\end{equation*}
Hence if $z\in B(R/2)$, we have
\begin{equation*}
\Psi(z)\leq C_{1}{1+R\over R^{2}}+C_{2}\left({f(\max_{B(R)}U)\over
R}\right)^{2}
\end{equation*}
Since $(\Omega,g)$ is complete, letting $R$ to infinity,  we obtain
$\Psi(z)=\Vert\nabla_{g} U(z)\Vert^{2}=0$ i.e $\det(u_{i\bar j})=C$
on $\Omega$.
\end{proof}

\end{document}